\documentclass[11pt]{amsart}
\usepackage{amssymb,latexsym,comment,url,amsmath}

\usepackage[T1]{fontenc}
\usepackage{rotating,graphicx}
\usepackage{currvita}
\usepackage{soul}
\usepackage{diagbox} 
\usepackage{xcolor}
\usepackage{subcaption}
\usepackage{multirow}
\usepackage{array}
\usepackage{longtable}
\usepackage{ltablex}
\usepackage{float}
\usepackage{lmodern,babel,adjustbox,booktabs,multirow}
\usepackage[T1]{fontenc}
\usepackage[utf8]{inputenc}
\usepackage{ upgreek }
\usepackage{caption}
\usepackage{placeins}
\restylefloat{table}
\usepackage{blindtext}

\usepackage{hyperref}
\usepackage{tabularray}
 \usepackage{lmodern}
    \usepackage{caption}
    \usepackage{booktabs}
    \usepackage{multirow}
    \usepackage{makecell}

\newcommand{\tors}{\operatorname{tors}}
\newcommand{\Tors}{\operatorname{Tors}}

\newcommand{\isom}{ \cong }

\newcommand{\Q}{{\mathbb Q}}

\newcommand{\Z}{{\mathbb Z}}

\newcommand{\Magma}{{\sf MAGMA}}

\newcommand{\Mathematica}{{\sf Mathematica}}

\newenvironment{Proof}{\par\noindent{\sc Proof:}}%
                      {\hspace*{\fill}\nobreak$\Box$\par\medskip}
                       {\hspace*{\fill}\nobreak$\Box$\par\medskip}

\newtheorem{Proposition}{Proposition}[section]
\newtheorem{Theorem}[Proposition]{Theorem}

\makeglossary
\newtheorem*{theorem*}{Theorem}

\theoremstyle{definition}

\newtheorem{Remark}[Proposition]{Remark}

\renewcommand{\baselinestretch}{1.1}

\begin{document}

\title[Torsion of Elliptic Curves
over Number Fields]%
{On Torsion Subgroups of Elliptic Curves
over Quartic, Quintic and Sextic Number Fields}

\author[M. Umut Kazancıoğlu]%
{Mustafa Umut Kazancıoğlu}
\address{Faculty of Engineering and Natural Sciences, Sabanc{\i} University, Tuzla, \.{I}stanbul, 34956 Turkey}
\email{mustafa.kazancioglu@sabanciuniv.edu}
\author[M. Sadek]%
{Mohammad~Sadek}
\email{mohammad.sadek@sabanciuniv.edu}

\begin{abstract}
The list of all groups that can appear as torsion subgroups of elliptic curves over number fields of degree $d$, $d=4,5,6$, is not completely determined. However, the list of 
 groups $\Phi^{\infty}(d)$, $d=4,5,6$, that can be realized as torsion subgroups for infinitely many non-isomorphic elliptic curves over these fields are known. We address the question of which torsion subgroups can arise over a given number field of degree $d$. In fact, given $G\in\Phi^{\infty}(d)$ and a number field $K$ of degree $d$, we give explicit criteria telling whether $G$ is realized finitely or infinitely often over $K$. We also give results on the field with the smallest absolute value of its discriminant such that there exists an elliptic curve with torsion $G$. Finally, we give examples of number fields $K$ of degree $d$, $d=4,5,6$, over which the Mordell-Weil rank of elliptic curves with prescribed torsion is bounded from above. 
 \end{abstract}
\maketitle

\let\thefootnote\relax\footnote{\textbf{Mathematics Subject Classification:} 11G05 \\

\textbf{Keywords:} torsion subgroups, modular curves, elliptic curves}

\section{Introduction}
Let $\Phi(d)$ be the set of possible isomorphism torsion structures $E(K)_{\tors}$, where $K$ runs
through all number fields $K$ of degree $d$ and $E$ runs through all elliptic curves defined over $K$.
 Let $\Phi^{\infty}(d)$ be the subset of isomorphic torsion structures in $\Phi(d)$ that occur infinitely often.
In other words, a torsion group $G \in \Phi^{\infty}(d)$ if there are infinitely many elliptic
curves $E$, non-isomorphic over $\overline{\Q}$, such that $E(K)_{\tors}\isom  G$ for some number field $K$ of degree $d$. We may identify elements of $\Phi(d)$ by pairs of positive integers $(m, mn)$, $m\ge 1$. Mazur, \cite{Maz78}, famously determined $\Phi(1)$; Kenku, Momose \cite{KM88} and Kamienny \cite{Kam92} determined $\Phi(2)$; whereas $\Phi(3)$ was recently determined in \cite{DEHMZB21}. More precisely, it is known that 
\begin{eqnarray*}
\Phi(1)&=&\{(1,n) : 1\leq n \leq 12, n\neq 11 \} \cup \{(2,2n) : 1\leq n \leq 4 \}, \\
\Phi(2)&=&\{(1,n) : 1\leq n \leq 18, n\neq 17 \} \cup \{(2,2n) : 1\leq n \leq 6 \} \cup \{(3,3), (3,6), (4,4)\},\\
\Phi(3)&=&\{(1,n) : 1\leq n \leq 21, n\neq 17,19 \} \cup \{(2,2n) : 1\leq n \leq 7 \}.
\end{eqnarray*}
Although $\Phi(d)$ has not been completely found for any $d\ge 4$, $\Phi^{\infty}(d)$ was given in \cite{Jeon} for $d=4$ and in \cite{DS16} when $ d=5, 6$, namely,
{\footnotesize\begin{eqnarray*}
    \Phi^{\infty}(4)&=&\{(1,n) : 1\leq n \leq 24, n\neq 19,23 \} \cup \{(2,2n) : 1\leq n \leq 9 \} 
\cup \{ (3,3n) : 1\leq n \leq 3  \} \cup \{(4,4), (4,8), (5,5), (6,6)\},\\
\Phi^{\infty}(5)&=&\{(1,n) : 1\leq n \leq 25, n\neq 23 \} \cup \{(2,2n) : 1\leq n \leq 8 \},\\
\Phi^{\infty}(6)&=&\{(1,n) : 1\leq n \leq 30, n\neq 23,25,29 \} \cup \{(2,2n) : 1\leq n \leq 10 \} 
\cup \{ (3,3n) : 1\leq n \leq 4  \} \cup \{(4,4), (4,8), (6,6)\}.
\end{eqnarray*}}
The latter result displays the groups that appear infinitely often as torsion subgroups of elliptic curves if one allows the base field to vary over all number fields of fixed degree $d$, $d=4,5,6$. Yet, fixing a number field $K$, $[K:\Q]=4,5$ or $6$, the aforementioned result cannot be used to determine the possible torsion subgroups of an elliptic curve defined over $K$. 

Although for any $(m,mn)\in\Phi(1)$ there are infinitely many non-isomorphic elliptic curves over $\Q$ whose torsion is $(m,mn)$, in particular $\Phi(d)=\Phi^{\infty}(d)$ for $d=1,2$. This is not generally the case over number fields of higher degree. In fact, F. Najman gave the first example showing that $\Phi(3)\neq\Phi^{\infty}(3)$ in \cite{Naj12}. More precisely, if $(m,mn)\in\Phi(d)$, $d\ge 3$. There are three possibilities for $(m,mn)$: 
\begin{itemize}
\item[i)] if $(m,mn)$ is the torsion of an elliptic curve over a number field $K$ of degree $d$, then there are 
infinitely many non-isomorphic elliptic curves over $K$ whose torsion is $(m,mn)$. In particular, for any number field $K$ of degree $d$, it is either that there are no elliptic curves over $K$ with torsion $(m,mn)$, or there are infinitely many non-isomorphic elliptic curves over $K$ with torsion $(m,mn)$,
\item[ii)] $(m,mn)$ is the torsion subgroup of at least one elliptic curve and at most finitely many elliptic curves over some number fields of degree $d$, while $(m,mn)$ is the torsion subgroup of infinitely
many non-isomorphic elliptic curves over the other number fields of degree $d$,
\item[iii)] $(m,mn)$ is the torsion subgroup for at most finitely many non-isomorphic
elliptic curves over each field of degree $d$. 
\end{itemize}
Given $(m,mn)\in \Phi^{\infty}(d)$, $m\ge 1$, $d=4,5,6$, we are interested in describing when one of the possibilities above is occurring for $(m,mn)$. This problem has been resolved for $d=2$ and $3$, see \cite{KamNaj,Naj11}. In fact, if $(m,mn)\not\in\Phi(1)$, then $(m,mn)$ can appear as the torsion subgroup
of infinitely many non-isomorphic elliptic curves over a number field $K$ only if the corresponding modular curve $X_1(m,mn)$ has
genus at most $1$. If the genus of $X_1(m,mn)$ is at least $2$, then $(m,mn)$ occurs as the torsion subgroup
of at most finitely many non-isomorphic elliptic curves over $K$.

If the genus of $X_1(m,mn)$ is one, we give a complete description of the number fields of degree $d$, $d=4,5,6$, over which the torsion subgroup of $X_1(m,mn)$ grows. We then classify these torsion points into cuspidal or non-cuspidal points and we display their fields of definition. If the rank of $X_1(m,mn)$ is zero over these number fields, then there 
are finitely many elliptic curves with torsion $(m,mn)$. In addition, Weierstrass equations describing these elliptic curves can be given explicitly. This, together with the study of genus zero modular curves, will enable us to decide which groups $(m,mn)\in\Phi^{\infty}(d)$ occur infinitely, respectively finitely, many often over a fixed number field $K$ of degree $d$, $d=4,5,6$



Over the rational field, there exists an elliptic curve with every possible torsion and 
rank $0$. Also, it is conjectured that one can find elliptic curves with arbitrarily large rank and any possible torsion over the rational field. A phenomenon that was observed in \cite{Dujella} and \cite{KamNaj} is that Mordell-Weil ranks of elliptic curves with prescribed torsion is bounded from above 
over certain quadratic, cubic and quartic fields. We substantiate this observation by giving examples of quartic, quintic and sextic number fields over which the Mordell-Weil rank of elliptic curves with certain torsion is bounded from above. 

 For each torsion group $(m,mn)$, such that $X_1(m,mn)$ has genus at most $1$, we find the number field of degree $d=4,5,6$ with the smallest absolute value of the discriminant such that an elliptic curve with torsion $(m,mn)$ occurs. We accomplish this by searching
through number fields with ascending discriminant until we hit a number field over which the
torsion $(m,mn)$ is possible.

All the necessary information needed to investigate the growth of the torsion of the elliptic modular curves $X_1(m, mn)$, which are discussed in this article, can be found on the LMFDB database webpage, \cite{lmfdb}. The study of the growth of the torsion subgroups of these curves are obtained from the analysis presented in \cite{GN20} and \cite{GN23}.

\subsection*{Acknowledgment}
The authors are indebted to the anonymous referee for the thorough reading of the manuscript and for many suggestions and comments that improved the manuscript. M.U. Kazancıo\u{g}lu would like to thank Andrew Sutherland for his help regarding calculations related to the curve $X_1(3,9)$. This work is based on the master's thesis of M. U. Kazanc{\i}o\u{g}lu, \cite{Kaz2}. All the calculations in this work were performed using \Magma, \cite{Magma}, and \Mathematica, \cite{Mathematica}. All codes used throughout the article can be found in \cite{Kaz}.

  This work is supported by The Scientific and Technological Research Council of Turkey, T\"{U}B\.{I}TAK, research grant ARDEB 1001/122F312. M. Sadek acknowledges the support of BAGEP Award of the Science Academy, Turkey.

\section{Preliminaries} 

We  recall that $K$-rational points of the affine modular curve $Y_1(m,mn)$ determine isomorphism classes of triples $(E,P_m,P_{mn})$, where $E$ is an elliptic curve over $K$, $P_m$ and $P_{mn}$ are generators of the subgroup of $E$ which is isomorphic to $(m,mn)$. 
In the case $m=1$, we write $Y_1(n)$ for $Y_1(1,n)$. 
The curve $X_1(m,mn)$ is a compactification of $Y_1(m,mn)$ obtained by adjoining its cusps.

The modular curve $X_1(m,mn)$ is defined over the cyclotomic field 
$\mathbb{Q}(\zeta_m)$, see \cite{DS16}. Despite that there might be points on the latter modular curve defined over smaller fields, these points do not give rise to elliptic curves with the desired torsion over these smaller fields. Explicit models for $X_1(m,mn)$, $m\ge 1$, can be found in \cite{Drew}. 

If $X_1(m,mn)$ is an elliptic curve, then we write $\Tors(X_1(m,mn)(K))$ for its $K$-rational torsion subgroup.

In the following tables, we list down the genus of the curve $X_1(m,mn)$, where $(m,mn)$ is either $(1,n)$, $4\le n\le 30$; or $(2,2n)$, $1\le n\le 10$; or $(3,3n)$, $1\le n\le 4$; or $(4,4n)$, $n=1, 2$; or $(m,m)$, $m=5,6$. We notice that if $(m,mn)$ is one of the groups appearing in $\Phi(d)$, $d\le 3$, or $\Phi^{\infty}(d)$, $4\le d\le 6$, then $(m,mn)$ will appear in these tables. The genus values can be found for example in \cite{Baaziz, Drew}.

{\scriptsize$$\begin{array}{|c||c|c|c|c|c|c|c|c|c|c|c|c|c|c|c|c|c|c|c|c|c|c|c|c|c|c|c|}
\hline
n & 4&5&6&7&8&9&10&11&12&13&14&15&16& 17&18&19&20&21&22&23&24&25&26&27&28&29&30\\
\hline
 g(X_1(n)) &0&0&0&0&0&0&0&1& 0&2&1&1&2&5&2&7&3&5&6&12&5& 12&10&13&10&22&9\\
\hline
\end{array}
$$}

{\tiny$$\begin{array}{|c||c|c|c|c|c|c|c|c|c|c|c|c|c|c|c|c|c|c|c|c|c|c|c|c|c|c|c|}
\hline
(m,mn)& (2,2)& (2,4)&(2,6)&(2,8)&(2,10)& (2,12)& (2,14)& (2,16)& (2,18) & (2,20)\\
\hline
g(X_1(m,mn))&0&0&0&0&1&1&4&5&7&9\\
\hline
\end{array}
$$}

{\tiny$$\begin{array}{|c||c|c|c|c|c|c|c|c|c|c|c|c|c|c|c|c|c|c|c|c|c|c|c|}
\hline
(m,mn)&(3,3)&(3,6)&(3,9) & (3,12) & (4,4)& (4,8)& (5,5) & (6,6)\\
\hline
g(X_1(m,mn))& 0&0&1&3&0&1&0&1\\
\hline
\end{array}
$$}

Now it is obvious that any $(m,mn)\in\Phi(1)$, i.e., a group that is realised as a torsion subgroup of an elliptic curve over $\Q$, also occurs infinitely often over every number field of any degree. The reason is that the modular curves $X_1(m,mn)$ are of genus $0$. For explicit construction of elliptic curves over a number field $K$ whose torsion subgroup lies in $\Phi(1)$, see \cite{Ku79}. We also notice that if the genus of $X_1(m,mn)$, $m\ge 1$, is at least $2$, then the group $(m,mn)$ will appear at most finitely many times over any fixed number field by Faltings' celebrated theorem, \cite{Fal83}.

\section{The modular curves $X_1(n)$, $n=11, 14, 15$, and $X_1(2,2m)$, $m=5,6$}

We remark that the modular curves $X_1(m,mn)$ appearing in the title of this section are all elliptic curves defined over $\Q$.  

\subsection{The strategy}
If $X_1(m,mn)$ is an elliptic curve, we write $L_{(m,mn),d}^{(s,st)}$ for a number field of degree $d$ over which $X_1(m,mn)$ possesses torsion subgroup $(s,st)$. In fact, for $d=4,5,6$, it will be seen that there are at most two such number fields; and we distinguish between them by writing $L_{(m,mn),d}^{(s,st)'}$ for the second such number field. If $m=1$, we shortly write $L_{n,d}^{(s,st)}$. For explicit description of these fields, see Table \ref{Tab1}.

 Let $E$ be an elliptic curve defined over a number field $K$. It is known that the torsion subgroup, $E(K)_{\tors}$, of $E(K)$ is of the form $\Z/m\Z \times \Z/mn\Z$ for some positive integers
 $m $, $n$. If $F$ is a finite field extension of $K$, then one obvious question to ask is how $E(K)_{\tors}$ grows to $E(F)_{\tors}$. If $E(K)_{\tors}\subsetneq E(F)_{\tors}$, $F$ is said to be a number field over which a growth of the torsion occurs. The interested reader may consult \cite{GN20, GN23} and the references there for recent results on torsion growth of elliptic curves defined over number fields. 

Let $(m,mn)$ be such that $X_1(m,mn)$ is an elliptic curve defined over $\Q$. In what follows we describe the general strategy we use to obtain our findings. 
\begin{itemize}
    \item[i)] We check the number fields of degree $d$, $d=4,5,6$, over which a growth of the torsion of $X_1(m,mn)$ occurs. These number fields can be found explicitly in \cite{lmfdb}. We listed these number fields in Table \ref{Tab1}.
    \item[ii)] For each number field $K$ found in i), we explicitly compute the points in $\Tors(X_1(m,mn),K)$ with minimal field of definition of degree $d$.
    \item[iii)] We compute Weierstrass equations for the elliptic curves with torsion $(m,mn)$ over $K$, up to isomorphism, corresponding to the non-cuspidal points in $\Tors(X_1(m,mn),K)$. The non-cuspidal points and the corresponding elliptic curve can be found in Tables \ref{Tab2}, \ref{Tab3}, \ref{Tab4}, \ref{Tab5}, \ref{Tab6}, \ref{Tab7}.
    \item[iv)] We compute the Mordell-Weil rank of $X_1(m,mn)$ over $K$ to investigate the possibility of existence of infinitely many elliptic curves over $K$ with torsion $(m,mn)$ if the rank is positive. 
\end{itemize}

\subsection{The elliptic curve $X_1(11)$}
The modular curve $X_1(11)$ is an elliptic curve defined by the Weierstrass equation $y^2-y=x^3-x^2$ where $\Tors(X_1(11),\Q)\simeq\Z/5\Z$. The number fields $K$ of degree $d$, $d=4$, $5$ or $6$, over which a growth of torsion of $X_1(11)$ occurs can be seen in Table \ref{Tab1}.
\begin{Theorem}
\label{thm11}
\begin{itemize}
\item[a)] Over any quartic number field there is either no elliptic curve with torsion $(1,11)$ or there are infinitely many such curves. 
\item[b)] There are exactly three non-isomorphic elliptic curves with torsion $(1,11)$ over $L_{11,5}^{(1,25)}$. Over any other quintic number field there
are either no elliptic curves with torsion $(1,11)$ or there are infinitely many such curves.
\item[c)] If $[K:\Q]=6$, then
\begin{itemize}
\item[i)] if $K\supset L_{11,3}^{(1,10)}$ and $K\not\simeq L_{11,6}^{(2,10)}$, there is exactly one elliptic curve defined over $L_{11,3}^{(1,10)}$ whose $L_{11,3}^{(1,10)}$-torsion is $(1,11)$ over $K$, in particular, over any such number field $K$, there is either exactly one elliptic curve, up to isomorphism, with torsion $G$ such that $\Z/11\Z\hookrightarrow G$ or there are infinitely many elliptic curves with torsion $(1,11)$;
\item[ii)] if $K\simeq L_{11,6}^{(2,10)}$, there are exactly three elliptic curves defined over $L_{11,3}^{(1,10)}$ with torsion $(1,11)$ over $L_{11,6}^{(2,10)}$ where these curves are non-isomorphic over $L_{11,6}^{(2,10)}$;
\item[iii)] over any other sextic number field, there
are either no elliptic curves with torsion $(1,11)$ or there are infinitely many such curves.
\end{itemize}
\end{itemize}
\end{Theorem}
\begin{Proof}
Statement a) follows from Table \ref{Tab1} since there is no quartic extension over which the torsion subgroup of $X_1(11)$ grows. As for b), one may see in Table \ref{Tab1} that the only quintic number field over which the torsion subgroup of $X_1(11)$ grows is $L_{11,5}^{(1,25)}$. In fact, the torsion grows from $(1,5)$ to $(1,25)$. One may also check that the Mordell-Weil rank of $X_1(11)$ over $L_{11,5}^{(1,25)}$ is zero. In addition, non-cuspidal points of $X_1(11)$ over the latter quintic number field generates three non-isomorphic elliptic curves with torsion $(1,11)$, see Table \ref{Tab2}. Since there is no growth in the torsion subgroup of $X_1(11)$ over any other quintic number field, statement b) follows.

 The $L_{11,3}^{(1,10)}$-Mordell Weil rank of $X_1(11)$ is zero. In addition, the non-cuspidal torsion points of $X_1(11)$ over $L_{11,3}^{(1,10)}$ give rise to one elliptic curve with torsion $\Z/11\Z$, see Table  \ref{Tab2}, hence c)-i) follows. 

Finally, the Mordell-Weil rank of $X_1(11)$ over $L_{11,6}^{(2,10)}$ is zero. The non-cuspidal $L_{11,3}^{(1,10)}$-rational points in $\Tors(X_1(11),L_{11,6}^{(2,10)} )$ generate three elliptic curves defined over $L_{11,3}^{(1,10)}$ which are non-isomorphic over $L_{11,6}^{(2,10)}$, see Table \ref{Tab2}. Since there is no growth in the torsion subgroup of $X_1(11)$ over sextic number fields other than $L_{11,6}^{(2,10)}$, this completes the proof.
\end{Proof}

\subsection{The elliptic curve $X_1(14)$} The modular curve $X_1(14)$ is defined by the Weierstrass equation $y^2+xy+y=x^3-x$ where $\Tors(X_1(14),\Q)\simeq\Z/6\Z$. The number fields $K$ of degree $d$, $d=4$, $5$ or $6$, over which a growth of torsion of $X_1(14)$ occurs can be seen in Table \ref{Tab1}.

\begin{Theorem}
\label{thm14}
\begin{itemize}
\item[a)] If $[K:\Q]=4$, then
\begin{itemize}
\item[i)] if $L_{14,2}^{(2,6)}\subset K$, then there are exactly two non-isomorphic elliptic curves with $L_{14,2}^{(2,6)}$-rational torsion $(1,14)$, in particular, over any such number field $K$, there are either exactly two non-isomorphic elliptic curves with torsion $G$ such that $\Z/14\Z\hookrightarrow G$ or there are infinitely many elliptic curves with torsion $(1,14)$;
\item[ii)] there are exactly two non-isomorphic elliptic curves over $L_{14,4}^{(1,12)}$ whose $L_{14,4}^{(1,12)}$-rational torsion subgroup is $(1,14)$;
\item[iii)] over any other quartic number field, there is either no elliptic curve with torsion $(1,14)$ or there are infinitely many such curves. 
\end{itemize}
\item[b)] Over any quintic number field there is either no elliptic curve with torsion $(1,14)$ or there are infinitely many such curves. 
\item[c)] If $[K:\Q]=6$, then
\begin{itemize}
\item[i)] if $L_{14,2}^{(2,6)}\subset K$ and $K\not\simeq L_{14,6}^{(2,18)}$, then there are exactly two non-isomorphic elliptic curves with $L_{14,2}^{(2,6)}$-rational torsion $(1,14)$, in particular, over any such number field $K$, there are either exactly two non-isomorphic elliptic curves with torsion $G$ such that $\Z/14\Z\hookrightarrow G$ or there are infinitely many elliptic curves with torsion $(1,14)$;
\item[ii)] if $L_{14,3}^{(1,18)}\subset K$ and $K\not\simeq L_{14,6}^{(2,18)}$, then there are exactly two non-isomorphic elliptic curves with $L_{14,3}^{(1,18)}$-rational torsion $(1,14)$, in particular, over any such number field $K$, there are either exactly two non-isomorphic elliptic curves with torsion $G$ such that $\Z/14\Z\hookrightarrow G$ or there are infinitely many elliptic curves with torsion $(1,14)$;
\item[iii)] if $K\simeq  L_{14,6}^{(3,6)}$, then there are infinitely many non-isomorphic elliptic curve over $K$ with torsion $(1,14)$;
\item[iv)] if $K\simeq  L_{14,6}^{(1,18)}$, then there are exactly three non-isomorphic elliptic curve over $K$ with torsion $(1,14)$; and exactly one elliptic curve with torsion $(2,14)$;
\item[v)] if $K\simeq  L_{14,6}^{(2,18)}$, then there are exactly five non-isomorphic elliptic curve over $K$ with torsion $(1,14)$; and exactly one elliptic curve with torsion $(2,14)$;
\item[vi)] over any other sextic number field, there is either no elliptic curve with torsion $(1,14)$ or there are infinitely many such curves. 
\end{itemize}
\end{itemize}
\end{Theorem}
\begin{Proof} The proof is similar to the proof of Theorem \ref{thm11} by using Table \ref{Tab1} and Table \ref{Tab3}.
\end{Proof}

\begin{Remark}
According to Theorem \ref{thm14}, one sees that $|X_1(2,14)(L_{14,6}^{(1,18)})|\geq 1$ and $|X_1(2,14)(L_{14,6}^{(2,18)})| \geq 1$. 
\end{Remark}
\subsection{The elliptic curve $X_1(15)$}
The modular curve $X_1(15)$ is defined by the Weierstrass equation $y^2+xy+y=x^3+x^2$ where $\Tors(X_1(15),\Q)\simeq\Z/4\Z$. The number fields $K$ of degree $d$, $d=4$, $5$ or $6$, over which a growth of torsion of $X_1(15)$ occurs can be seen in Table \ref{Tab1}.

\begin{Theorem}
\label{thm15}
\begin{itemize}
\item[a)] If $[K:\Q]=4$, then
\begin{itemize}
\item[i)] if $K$ is isomorphic to either of the number fields $ L_{15,4}^{(1,16)}$, $L_{15,4}^{(1,16)'}$ or $L_{15,4}^{(2,8)}$, then there are exactly two non-isomorphic elliptic curves over $K$ with $K$-rational torsion $(1,15)$;
\item[ii)] if $L_{15,2}^{(2,4)}\subset K$ and $K\not\simeq L_{15,4}^{(2,8)}$, then there is exactly one elliptic curve, up to isomorphism, with $L_{15,2}^{(2,4)}$-rational torsion $(1,15)$, in particular, over any such number field $K$, there is either exactly one elliptic curve, up to isomorphism, with torsion $G$ such that $\Z/15\Z\hookrightarrow G$ or there are infinitely many elliptic curves with torsion $(1,15)$;
\item[iii)] if $L_{15,2}^{(1,8)}\subset K$ and $K$ is not isomorphic to any of the number fields $ L_{15,4}^{(1,16)}$, $L_{15,4}^{(1,16)'}$ or $L_{15,4}^{(2,8)}$, then there is exactly one elliptic curve with $L_{15,2}^{(1,8)}$-rational torsion $(1,15)$, in particular, over any such number field $K$, there is either exactly one elliptic curve, up to isomorphism, with torsion $G$ such that $\Z/15\Z\hookrightarrow G$ or there are infinitely many elliptic curves with torsion $(1,15)$;
\item[iv)] over any other quartic number field, there is either no elliptic curve with torsion $(1,15)$ or there are infinitely many such curves. 
\end{itemize}
\item[b)] Over any quintic number field there is either no elliptic curve with torsion $(1,15)$ or there are infinitely many such curves. 
\item[c)] If $[K:\Q]=6$, then
\begin{itemize}
\item[i)] if $L_{15,2}^{(2,4)}\subset K$, then there is exactly one elliptic curve with $L_{15,2}^{(2,4)}$-rational torsion $(1,15)$, in particular, over any such number field $K$, there is either exactly one elliptic curve, up to isomorphism, with torsion $G$ such that $\Z/15\Z\hookrightarrow G$ or there are infinitely many elliptic curves with torsion $(1,15)$;
\item[ii)] if $L_{15,2}^{(1,8)}\subset K$, then there is exactly one elliptic curve with $L_{15,2}^{(1,8)}$-rational torsion $(1,15)$, in particular, over any such number field $K$, there is either exactly one elliptic curve, up to isomorphism, with torsion $G$ such that $\Z/15\Z\hookrightarrow G$ or there are infinitely many elliptic curves with torsion $(1,15)$;
\item[iii)] over any other sextic number field, there is either no elliptic curve with torsion $(1,15)$ or there are infinitely many such curves. 
\end{itemize}
\end{itemize}
\end{Theorem}
\begin{Proof} The proof is analogous to the proof of Theorem \ref{thm11}, see Table \ref{Tab1} and Table \ref{Tab4}.
\end{Proof}


\subsection{The elliptic curve $X_1(2,10)$}
A Weierstrass equation describing the elliptic curve $X_1(2,10)$ is given by $y^2=x^3+x^2-x=x(x^2+x-1).$ In addition, $\Tors(X_1(2,10),\mathbb{Q}) \simeq \mathbb{Z}/ 6\mathbb{Z}.$  The number fields $K$ of degree $d$, $d=4$, $5$ or $6$, over which a growth of torsion of $X_1(2,10)$ occurs can be seen in Table \ref{Tab1}.

\begin{Theorem}
\label{thm210}
\begin{itemize}
\item[a)] If $[K:\Q]=4$, then
\begin{itemize}
\item[i)] if $K\simeq L_{(2,10),4}^{(1,12)}$, then there is exactly one elliptic curve over $K$, up to isomorphism, with $K$-rational torsion $(2,10)$;
\item[ii)] over any other quartic number field, there is either no elliptic curve with torsion $(2,10)$ or there are infinitely many such curves. 
\end{itemize}
\item[b)] Over any quintic number field there is either no elliptic curve with torsion $(2,10)$ or there are infinitely many such curves. 
\item[c)] If $[K:\Q]=6$, then
\begin{itemize}
\item[i)] if $K\simeq L_{(2,10),6}^{(3,6)}$, then there is exactly one elliptic curve, up to isomorphism, with $K$-rational torsion $(2,10)$;
\item[ii)] over any other sextic number field, there is either no elliptic curve with torsion $(2,10)$ or there are infinitely many such curves. 
\end{itemize}
\end{itemize}
\end{Theorem}

\begin{Proof}
The proof is similar to the proof of Theorem \ref{thm11}, see Table \ref{Tab1} and \ref{Tab5}.
\end{Proof}

\subsection{The elliptic curve $X_1(2,12)$}

A Weierstrass equation describing the elliptic curve $X_1(2,12)$ is defined over $\Q(\sqrt{-1})$ by $y^2=x(x^2-x+1).$ Further, $\Tors(X_1(2,12),\mathbb{Q}) \simeq \mathbb{Z}/ 4\mathbb{Z}.$ The number fields $K$ of degree $d$, $d=4$, $5$ or $6$, over which a growth of torsion of $X_1(2,12)$ occurs can be seen in Table \ref{Tab1}.

\begin{Theorem}
\label{thm212}
If $[K:\Q]=4,5$ or $6$, then 
 there is either no elliptic curve over $K$ with $K$-torsion $(2,12)$ or there are infinitely many such curves. 
\end{Theorem}
\begin{Proof}
The proof is similar to the proof of Theorem \ref{thm11}, see Table \ref{Tab1}.
\end{Proof}

\section{The modular curves $X_1(3,9)$, $X_1(4,8)$ and $X_1(6,6)$}
We recall that the modular curve $X_1(m,mn)$ is defined over the cyclotomic field $\Q(\zeta_m)$.  This holds because if the base field $K$ over which an elliptic curve $E$ is defined does not contain a primitive $m$-th root of unity, then $(m,mn)$ can not be contained in the $K$-torsion subgroup of $E$. It follows that the curves $X_1(3,9)$, $X_1(4,8)$ and $X_1(6,6)$ are defined over $\Q(\zeta_3)$, $\Q(\sqrt{-1})$ and $\Q(\zeta_3)$, respectively. 
\subsection{The elliptic curve $X_1(3,9)$}
The elliptic curve $X_1(3,9)$ is defined by the Weierstrass equation $y^2+y=x^3$ where $\Tors(X_1(3,9),\Q)\simeq\Z/3\Z$. The number fields $K$ of degree $d$, $d=4$, $5$ or $6$, over which a growth of torsion of $X_1(3,9)$ occurs can be seen in Table \ref{Tab1}. 
\begin{Theorem}
\label{thm39}
\begin{itemize}
\item[a)] If $[K:\Q]=4$ or $5$, then 
 there is either no elliptic curve over $K$ with $K$-torsion $(3,9)$ or there are infinitely many such curves. 
 \item[b)] If $[K:\Q]=6$, then 
 \begin{itemize}
 \item[i)] if $K\simeq L_{(3,9),6}^{(3,9)}$ or $L_{(3,9),6}^{(6,6)}$, then there is exactly one elliptic curve, up to isomorphism, over $K$ with torsion $(3,9)$,
 \item[ii)] over any other sextic number field, there is either no elliptic curve with torsion $(3,9)$ or there are infinitely many such curves.
 \end{itemize}
\end{itemize}
\end{Theorem}
\begin{Proof}
The proof is similar to the proof of Theorem \ref{thm11}, see Table \ref{Tab1} and Table \ref{Tab6}.

\end{Proof}

\subsection{The elliptic curve $X_1(4,8)$}
The elliptic curve $X_1(4,8)$ is defined by the Weiertrass equation $y^2=x^3-x$ over $\Q(\sqrt{-1})$. 
We have that
$\Tors(X_1(4,8),\mathbb{Q})\simeq\mathbb{Z} / 2\mathbb{Z} \times  \mathbb{Z} /2\mathbb{Z}.$  The number fields $K$ of degree $d$, $d=4$, $5$ or $6$, over which a growth of torsion of $X_1(4,8)$ occurs, can be seen in Table \ref{Tab1}.

\begin{Theorem}
\label{thm212}
If $[K:\Q]=4,5$ or $6$, then 
 there is either no elliptic curve over $K$ with $K$-torsion $(4,8)$ or there are infinitely many such curves. 
\end{Theorem}
\begin{Proof}
 The proof is similar to the proof of Theorem \ref{thm11}, see Table \ref{Tab1}.
\end{Proof}

\subsection{The elliptic curve $X_1(6,6)$}
The curve $X_1(6,6)$ is defined by the equation $y^2=x^3+1$ where $\Tors(X_1(6,6),\Q)\simeq\Z/6\Z$. The number fields $K$ of degree $d$, $d=4$, $5$ or $6$, over which a growth of torsion of $X_1(6,6)$ occurs can be seen in Table \ref{Tab1}.

\begin{Theorem}
\label{thm66}
\begin{itemize}
\item[a)] If $[K:\Q]=4$ or $5$, then 
 there is either no elliptic curve over $K$ with $K$-torsion $(6,6)$ or there are infinitely many such curves. 
 \item[b)] If $[K:\Q]=6$, then 
 \begin{itemize}
 \item[i)] if $K\simeq L_{(6,6),6}^{(6,6)}$, then there is exactly one elliptic curve, up to isomorphism, over $K$ with torsion $(6,6)$,
 \item[ii)] over any other sextic number field, there is either no elliptic curve with torsion $(6,6)$ or there are infinitely many such curves.
 \end{itemize}
\end{itemize}
\end{Theorem}
\begin{Proof}
The proof is analogous to the proof of Theorem \ref{thm11}, see Table \ref{Tab1} and Table \ref{Tab7}.

\end{Proof}
\section{The modular curves $X_1(3,3)$, $X_1(3,6)$, $X_1(4,4)$ and $X_1(5,5)$}
The modular curves $X_1(3,3)$, $X_1(3,6)$, $X_1(4,4)$ and $X_1(5,5)$ are curves of genus $0$.
\subsection{The genus $0$ curve $X_1(3,3)$}
Recall that $X_1(3,3)$ is defined over the field $\mathbb{Q}(\zeta_3)$  where $\zeta_n$ is an $n$-th primitive root of unity, $n\ge 3$.
A Weierstrass equation of an elliptic curve with torsion subgroup $(3,3)$, \cite{Drew}, is the following:
$$\mathcal{E}_v(3,3): y^2+ ((\zeta_3+2)v+(1-\zeta_3))xy + ((\zeta_3+1)v^2-\zeta_3v)y=x^3,\qquad v \in K.$$
The discriminant of the latter curve is given by $$\Updelta(3,3)=3 (v-1)^2 v^3 \left((-1)^{2/3} v+v-(-1)^{2/3}\right)^3 \left(-2 \sqrt[3]{-1} v+4 (-1)^{2/3} v+3 v-\sqrt[3]{-1}-(-1)^{2/3}\right).$$
Notice that $\Updelta(3,3)=0$ for finitely many $v$-values, namely, $v=0,1,\sqrt[3]{-1}$. For any other $v$-value, the torsion subgroup of the elliptic curve $\mathcal{E}_v(3,3)$ is $(3,3)$.
\subsection{The genus $0$ curve $X_1(3,6)$}
The curve $X_1(3,6)$ is defined over $\mathbb{Q}(\zeta_3)$.
A Weierstrass equation of an elliptic curve with torsion subgroup $(3,6)$, \cite{Drew}, is the following:
$$\mathcal{E}_v(3,6): y^2+ (t+2)xy + (-t(t+1))y=x^3+(-t(t+1))x^2,\qquad t=\frac{4v^2+6v+3}{v^3},\, v \in K.$$ 
The discriminant is given by
\begin{equation*} 
\Updelta(3,6)=\frac{(v+1)^4 (2 v+1) \left(v^2+3 v+3\right)^4 \left(4 v^2+3\right) \left(4 v^2+6 v+3\right)^3 \left(v^3-4 v^2-6 v-3\right)^2}{v^{30}}
\end{equation*}
Notice that $\Updelta(3,6)=0$ if and only if $v$ is one of the values $-1$, $-\frac{1}{2}$, $\pm\frac{ i \sqrt{3}}{2}$, 
$\frac{1}{4} \left(-3 \pm i \sqrt{3}\right)$, 
$\frac{1}{2} \left(-3\pm i \sqrt{3}\right)$, 
and a root of the polynomial $v^3-4 v^2-6 v-3$. For any other $v$-value, the corresponding elliptic curve has torsion subgroup $(3,6)$.
\subsection{The genus $0$ curve $X_1(4,4)$}
The curve $X_1(4,4)$ is defined over the field $\mathbb{Q}(\zeta_4)$.
An equation describing an elliptic curve with torsion subgroup $(4,4)$, \cite{Drew}, is given by:
$$\mathcal{E}_v(4,4): y^2+ xy -ty=x^3  -tx^2, \qquad t=\frac{(1-v)(v^2-2v+2)}{2v^4},\, v \in K.$$ 
The discriminant is given by $$\Updelta(4,4)=-\frac{27 (1-v)^4 \left(v^2-2 v+2\right)^4}{16 v^{16}}.$$
Notice that $\Updelta(4,4)=0$ if and only if $v=1, 1-i$ or $1+i$. For any other $v$-value, the corresponding elliptic curve has torsion subgroup $(4,4)$.
\subsection{The genus $0$ curve $X_1(5,5)$}
The curve $X_1(5,5)$ is defined over the field $\mathbb{Q}(\zeta_5)$. 
A Weierstrass equation describing an elliptic curve with torsion subgroup $(5,5)$, \cite{Drew}, is:
$$\mathcal{E}_v(5,5): y^2+ (1-t)xy -ty=x^3 -tx^2,\qquad t=\frac{U}{V(U+1)},$$ where
{\footnotesize $$U=\frac{(2-a)v^2 + (2-a)v + a+3}{5(v+1)}, \quad V=\frac{-((a+2)v^2 + (5a+9)v + (25a+41))}{(v^3 + (-3a-2)v^2 + (2a+6)v +(5a+9))},\qquad a=\frac{\zeta_5+1}{\zeta_5}.$$}
In addition, the vanishing of the discriminant $\Updelta(5,5)$ is equivalent to the vanishing of the following polynomial in $\Q[v]$.
\begin{equation*}
\begin{split}
    &(v+1)^4 \left(\gamma^{2} v^2-3 \left(1+\gamma^{2}+\gamma^{4}\right) v+9 \gamma^{2}+5 \gamma^{4}+5\right)^4 \\ & \left(-v^2 \left(-1+\gamma^{2}\right)^2-v \left(-1+\gamma^{2}\right)^2+3 \gamma^{2}+\gamma^{4}+1\right)^4f_{5}(v)f_{10}(v)
\end{split}
\end{equation*}
where $\gamma\ne -1$ is a fifth root of $-1$ and $f_{i}(v)\in\Q[v]$ is an irreducible polynomial of degree $i$, $i=5,10$.


 Notice that $\Updelta(5,5)=0$ if and only if $v=-1$, $1+\gamma+2 \gamma^{2}$, $2-\gamma+\gamma^{2}-3 \gamma^{3}$, $-\gamma-\gamma^{3}$, $-1+\gamma+\gamma^{3}$, or $v$ is root of $f_{i}(v)$, $i=5,10$. For any other $v$-value, the corresponding elliptic curve has torsion subgroup $(5,5)$.

\section{Rank of elliptic curves with a given torsion}
It has been believed that an elliptic curve with prescribed torsion over the rational field
can have arbitrary large rank. In fact, for any torsion group of an elliptic curve over the rational field, there exists an elliptic curve with that torsion group and
rank at least $3$, see \cite{Dujellawebpage}. Unlike the rational field, it appears that upper bounds can be found for the Mordell-Weil rank of elliptic curves with prescribed torsion over number fields of higher degree, see \cite{Dujella, KamNaj}.

In Theorems \ref{thm11}, \ref{thm14}, \ref{thm15}, \ref{thm210}, \ref{thm39}, \ref{thm66}, each of the elliptic modular curves $X_1(11)$, $X_1(14)$, $X_1(15)$, $X_1(2,10)$, $X_1(3,9)$, $X_1(6,6)$ turn out to have non-cuspidal torsion points over certain number fields of degree $d$, $d=2,3,4,5,6$. Moreover, the Mordell-Weil group of these modular curves have rank zero over these number fields. It follows that there are only finitely many elliptic curves with the corresponding torsion over these number fields. We use \Magma\, to compute the Mordell-Weil rank for each of these finitely many elliptic curves yielding an upper bound on the rank. 
\begin{Theorem}
\begin{itemize}
\item[a)] All elliptic curve with torsion $(1,11)$ over $L_{11,5}^{(1,25)}$ have either rank $0$ or $1$. 
 \item[b)] All elliptic curve with torsion $(1,11)$ over either fields $L_{11,3}^{(1,10)}$ or $L_{11,6}^{(2,10)}$ have rank $0$. 
 \item[c)] All elliptic curve with torsion $(1,14)$ over either fields $L_{14,2}^{(2,6)}$, $L_{14,3}^{(1,18)}$, $L_{14,4}^{(1,12)}$, $L_{14,6}^{(1,18)}$ or $L_{14,6}^{(2,18)}$ have rank $0$. 
 \item[d)] All elliptic curve with torsion $(1,15)$ over either fields $L_{15,2}^{(1,8)}$, $L_{15,2}^{(2,4)}$, $L_{15,4}^{(2,8)}$, $L_{15,4}^{(1,16)}$ or $L_{15,4}^{(1,16)'}$ have rank $0$. 
 \item[e)] All elliptic curve with torsion $(2,10) $ over either fields $L_{(2,10),4}^{(1,12)}$ or $L_{(2,10),6}^{(3,6)}$ have rank $0$.
 \item[f)] All elliptic curve with torsion $(3,9) $ over either fields $L_{(3,9),6}^{(6,6)}$ or $L_{(3,9),6}^{(3,9)}$ have rank $0$.
 \item[g)] All elliptic curve with torsion $(6,6) $ over $L_{(6,6),6}^{(6,6)}$ have rank $0$.
\end{itemize}
\end{Theorem}
\begin{Proof}
We recall that the Mordell-Weil rank of $X_1(11)$ over $L_{11,5}^{(1,25)}$ is zero. It follows that the only elliptic curves with torsion subgroup $(1,11)$ over $L_{11,5}^{(1,25)}$ correspond to the $L_{11,5}^{(1,25)}$-rational torsion points of $X_1(11)$. 
In Table \ref{Tab1}, one sees that the non-cuspidal torsion points of $X_1(11)$ over $L_{11,5}^{(1,25)}$ give rise to three non-isomorphic elliptic curves over $L_{11,5}^{(1,25)}$. In Table \ref{Tab2}, it can be seen that one of these curves has Mordell-Weil rank $1$ over $L_{11,5}^{(1,25)}$, whereas the other two elliptic curves have Mordell-Weil rank $0$ over $L_{11,5}^{(1,25)}$, hence statement a) follows. 

The proof of the other statements can be obtained in a similar fashion by considering the information provided in Table \ref{Tab1} and the last column of Tables \ref{Tab2},\ref{Tab3},\ref{Tab4},\ref{Tab5},\ref{Tab6}, \ref{Tab7}.
\end{Proof}

\section*{Appendix}
In Tables \ref{Tab9}, \ref{Tab8}, we list the minimal polynomial of $a\in \overline{\mathbb{Q}}$ such that there exists an elliptic curve over $\mathbb{Q}(a)$ with torsion subgroup $(m,mn)$ and the discriminant of $\mathbb{Q}(a)$ is minimal.
In fact, Given $(m,mn)\in \Phi^{\infty}(d)$, $d=4,5,6$, such that the genus of $X_1(m,mn)$ is at most one, we find the number field $K$ of degree $d$ with the smallest absolute value of the discriminant $|\Delta|$ over which $(m,mn)$ is realized as the torsion subgroup of an elliptic curve. We accomplish this by searching through number fields with ascending discriminant, see \cite{JR14}, until
we find a number field over which $X_1(m,mn)$ has a non-cuspidal point. 

In Tables \ref{Tab2},\ref{Tab3},\ref{Tab4},\ref{Tab5},\ref{Tab6} and \ref{Tab7}, we list the equations of the non-isomorphic elliptic curves over $\mathbb{Q}(a)$ with torsion subgroup $(1,11)$, $(1,14)$, $(1,15)$, $(2,10)$, $(3,9)$ and $(6,6)$, respectively.

Table \ref{Tab1} provides information related to the growth of the torsion subgroup of $X_1(m, mn)$, detailing the description of the growth over the specific number fields and the number of non-isomorphic elliptic curves with torsion subgroup $(m,mn)$ corresponding to the non-cuspidal torsion points of $X_1(m,mn)$.

\tiny{\begin{table}[H]
\caption{Minimal Discriminant Genus $1$}
\begin{center}
\begin{itemize}
    \item $1^{st}$ column: Modular elliptic curve $X_1(m,mn)$.
    \item $2^{nd}$ column: Torsion Subgroup of $X_1(m,mn)$ over $\mathbb{Q}$.
    \item $3^{rd}$ column: Torsion Subgroup of $X_1(m,mn)$ over $\mathbb{Q}(a)$.
    \item $4^{th}$ column: Discriminant of $\mathbb{Q}(a)$.
    \item $5^{th}$ column: Minimal polynomial for $a\in \overline{\mathbb{Q}}$ such that there exists an elliptic curve over $\mathbb{Q}(a)$ with torsion subgroup $(m,mn)$ and discriminant of $\mathbb{Q}(a)$ minimal.
    
\end{itemize}
\begin{tabular}{| m{2.3cm} | m{2.9cm} | m{3.2cm} | m{1.9cm} |  m{5cm} |  } 
  \hline 
  $X_1(11)$&$(1,5)$&$(1,5)$&$117$&$x^4-x^3-x^2+x+1$ \\ \hline
  $X_1(14)$&$(1,6)$&$(1,6)$&$144$&$x^4-x^2+1$ \\ \hline
  $X_1(15)$&$(1,4)$&$(1,16)$&$125$&$x^4-x^3+x^2-x+1$ \\ \hline
  $X_1(2,10)$&$(1,6)$&$(1,6)$&$117$&$x^4-x^3-x^2+x+1$ \\ \hline
  $X_1(2,12)$&$(1,4)$&$(2,4)$&$189$&$x^4-x^3+2x+1$ \\ \hline
  $X_1(3,9)$&$(1,3)$&$(3,3)$&$900$&$x^4-x^3+2x^2+x+1$ \\ \hline
  $X_1(4,8)$&$(2,2)$&$(2,4)$&$400$&$x^4+3x^2+1$ \\ \hline
  $X_1(6,6)$&$(1,6)$&$(2,6)$&$1764$&$x^4-x^3-x^2-2x+4$ \\ \hline \hline
  $X_1(11)$&$(1,5)$&$(1,5)$&$1649$&$x^5-x^4+x^2-x+1$ \\ \hline
  $X_1(14)$&$(1,6)$&$(1,6)$&$1609$&$x^5-x^3-x^2+x+1$ \\ \hline
  $X_1(15)$&$(1,4)$&$(1,16)$&$1609$&$x^5-x^3-x^2+x+1$ \\ \hline
  $X_1(2,10)$&$(1,6)$&$(1,6)$&$1777$&$x^5-x^4+x^3-2x^2+x+1$ \\ \hline
  $X_1(2,12)$&$(1,4)$&$(1,4)$&$1649$&$x^5-x^4+x^2-x+1$ \\ \hline\hline
  $X_1(11)$&$(1,5)$&$(1,5)$&$-10051$&$x^6-x^5+2x^4-2x^3+2x^2-2x+1$ \\ \hline
  $X_1(14)$&$(1,6)$&$(1,6)$&$-10571$&$x^6-x^5-x^4+2x^3-x+1$ \\ \hline
  $X_1(15)$&$(1,4)$&$(1,4)$&$-10051$&$x^6-x^5+2x^4-2x^3+2x^2-2x+1$ \\ \hline
  $X_1(2,10)$&$(1,6)$&$(1,6)$&$-10051$&$x^6-x^5+2x^4-2x^3+2x^2-2x+1$ \\ \hline
  $X_1(2,12)$&$(1,4)$&$(2,4)$&$-9747$&$x^6-x^5+x^4-2x^3+4x^2-3x+1$ \\ \hline
  $X_1(3,9)$&$(1,3)$&$(3,3)$&$-14283$&$x^6-x^5+x^4-2x^3+x^2+1$ \\ \hline
  $X_1(4,8)$&$(2,2)$&$(2,4)$&$-33856$&$x^6+x^4+2x^2+1$ \\ \hline
  $X_1(6,6)$&$(1,6)$&$(2,6)$&$-14283$&$x^6-x^5+x^4-2x^3+x^2+1$ \\ \hline

\end{tabular}
\begin{itemize}
    \item[] Note: For the modular curve $X_1(6,6)$, we have not been able to compute its Mordell-Weil rank over the number field $\mathbb{Q}(a)$ where the minimal polynomial of $a$ is $x^6-x^5+x^4-2x^3+4x^2-3x+1$. The latter number field has discriminant $-9747$. Using the \Magma\, command AnalyticRank(), the analytic rank of $X_1(6,6)$ over the number field $\mathbb{Q}(a)$ is $0$. Therefore, we may claim that the Mordell-Weil rank of $X_1(6,6)$ is conjecturally $0$ over this sextic number field. 
\end{itemize}
\label{Tab9}
\end{center}
\end{table}}

\tiny{\begin{table}[H]
\caption{Minimal Discriminant Genus $0$}
\begin{center}
\begin{itemize}
    \item $1^{st}$ column: Modular curve $X_1(m,mn)$.
    \item $2^{nd}$ column: Discriminant of $\mathbb{Q}(a)$.
    \item $3^{rd}$ column: Minimal polynomial for $a\in \overline{\mathbb{Q}}$ such that there exists an elliptic curve over $\mathbb{Q}(a)$ with torsion subgroup $(m,mn)$ and discriminant of $\mathbb{Q}(a)$ minimal.
    
\end{itemize}
\begin{tabular}{| m{4.2cm}  | m{1.9cm} |  m{5cm} |  } 
  \hline 
  $X_1(3,3)$ or $X_1(3,6)$&$117$&$x^4-x^3-x^2+x+1$ \\ \hline
  $X_1(4,4)$&$144$&$x^4-x^2+1$ \\ \hline
  $X_1(5,5)$&$125$&$x^4-x^3+x^2-x+1$ \\ \hline \hline
  $X_1(3,3)$ or $X_1(3,6)$&$-9747$&$x^6-x^5+x^4-2x^3+4x^2-3x+1$ \\ \hline
  $X_1(4,4)$&$-10816$&$x^6-x^4-2x^3+2x+1$ \\ \hline

\end{tabular}
\begin{itemize}
    \item[]Note: For the genus $0$ curves $X_(3,3)$, $X_1(3,6)$, $X_1(4,4)$ and $X_1(5,5)$ we have parametrization over $\mathbb{Q}(\zeta_3)$, $\mathbb{Q}(\zeta_3)$, $\mathbb{Q}(\zeta_4)$ and $\mathbb{Q}(\zeta_5)$, respectively. The above table shows the minimal number fields where those condition are satisfied. Note that $[\mathbb{Q}(\zeta_3):\mathbb{Q}]=[\mathbb{Q}(\zeta_4):\mathbb{Q}]=2$ and $[\mathbb{Q}(\zeta_5):\mathbb{Q}]=4$.
\end{itemize}
\label{Tab8}
\end{center}
\end{table}}

\tiny{\begin{table}[H]
   \caption{Equations of the non-isomorphic elliptic curves over $\mathbb{Q}(a)$ with torsion subgroup $(1,11)$.}
\begin{center}

\begin{tabular}{| m{5cm} |  m{12cm} | m{0.7cm} | } 
  \hline

 \multicolumn{3}{|l|}{Number field generated by a root of the polynomial $x^3-x^2+x+1$}\\
  \hline
  Point from $X_1(11)$ &Corresponding Elliptic Curve $E$ with torsion $\mathbb{Z} \slash 11\mathbb{Z}$ & Rank of $E$  \\ 
  \hline
$(-a , a^2 - a)$ & $y^2 + a^2xy + (4a^2 - 7a - 5)y = x^3 + (-5a^2 +a + 2)x^2$ & 0 \\ \hline \hline
\multicolumn{3}{|l|}{Number field generated by a root of the polynomial $x^5-x^4-4x^3+3x^2+3x-1$} \\ \hline
  Point from $X_1(11)$ &Corresponding Elliptic Curve $E$ with torsion $\mathbb{Z} \slash 11\mathbb{Z}$ & Rank of $E$  \\ 
  \hline
$(-a^3 + a^2 + 2a , -a^4 + a^3 + 2a^2)$ & $y^2 + (-4a^4 + 11a^3 - 3a^2 - 8a + 3)xy +
    (-652a^4 + 1739a^3 - 321a^2 - 1380a + 383)y = x^3 + (-155a^4 + 411a^3
    - 73a^2 - 325a + 90)x^2 $ & 0 \\ \hline
$(-3a^4 + 5a^3 +8a^2 - 14a + 3 , -11a^4 +20a^3 + 27a^2 - 55a + 14)$ & $y^2 + (501a^4 - 918a^3 -1241a^2 + 2537a -605)xy + (-7687269a^4 +14074083a^3 + 19055905a^2 - 38894025a +9252517)y = x^3 + (-484791a^4 + 887570a^3 + 1201744a^2 - 2452818a +583502)x^2$ & 0 \\ \hline
$(a^3 - a^2 - 3a +3 , -2a^4 + 3a^3 + 6a^2 - 8a)$ & $y^2 + (4a^4 - 8a^3 - 8a^2 + 22a - 10)xy +(-739a^4 + 1329a^3 + 1869a^2 - 3669a + 788)y = x^3 + (-180a^4 +318a^3 + 464a^2 - 877a +168)x^2$& 1\\ \hline \hline
\multicolumn{3}{|l|}{Number field generated by a root of the polynomial $x^6-x^5+2x^4-3x^3+2x^2-x+1$} \\ \hline
 
  Point from $X_1(11)$ &Corresponding Elliptic Curve $E$ with torsion $\mathbb{Z} \slash 11\mathbb{Z}$ & Rank of $E$  \\ 
  \hline
$(1/2(-a^5 + a^4 - a^3 + 2a^2) , 1/2)$ & $y^2 + 1/4(-3a^5 + 3a^4 - 3a^3 + 6a^2 - 1)xy +
1/32(a^4 + a^2 - 1)y = x^3 + 1/16(a^5 + a^3 - a^2 + 1)x^2$ & 0 \\ \hline
$(a^3 - a^2 + a , -a^5+ a^4 - a^3 + a^2 - a + 1)$ & $y^2 + (2a^5 - 2a^4 + 3a^3 - 3a^2 + 3a)xy +(-34a^5 + 34a^4 - 6a^3 + 6a^2 - 6a - 10)y = x^3 + (-12a^5 + 12a^4 -2a^3 + 2a^2 - 2a -4)x^2$ & 0 \\ \hline
$(-a^5 + a^4 - 2a^3 + 3a^2 - a + 2 , -a^5 -a^3 + 2a^2 + a + 2)$ & $y^2 + (a^5 + a^4 - a^2 - 3a)xy + (-62a^5 + 28a^4- 90a^3 + 152a^2 + 6a + 80)y = x^3 + (-22a^5 + 10a^4 - 32a^3 + 54a^2+ 2a +28)x^2$ & 0 \\ \hline 
\end{tabular}
  \label{Tab2}
\end{center}
\end{table}}

\tiny{\begin{table}[H]
   \caption{Equations of the non-isomorphic elliptic curves over $\mathbb{Q}(a)$ with torsion subgroup $(1,14)$.}
\begin{center}

{\tiny 
\begin{tabular}{| m{5cm} |  m{13cm} |  m{0.6cm} |  } 
  \hline
  \multicolumn{3}{|l|}{Number field generated by a root of the polynomial $x^2-x+2$} \\ \hline
   Point from $X_1(14)$&Corresponding Elliptic Curve $E$ with torsion $\mathbb{Z} \slash 14\mathbb{Z}$ & Rank of $E$ \\ 
  \hline
$(1/4(-a + 2) , 1/8(a - 6))$ & $y^2 + 1/28(a + 31)xy + 1/56(a + 5)y = x^3 +1/56(a +5)x^2$ & 0\\ \hline
$(-a,-2)$ & $y^2 + 1/7(-4a + 9)xy + 1/7(-2a + 2)y = x^3 +1/7(-2a + 2)x^2$& 0 \\ \hline \hline
\multicolumn{3}{|l|}{Number field generated by a root of the polynomial $x^3-x^2-2x+1$} \\ \hline

  Point from $X_1(14)$ &Corresponding Elliptic Curve $E$ with torsion $\mathbb{Z} \slash 14\mathbb{Z}$ & Rank of $E$ \\ 
  \hline
$(2a^2 + 2a - 1 , 6a^2 + 4a - 4)$ & $y^2 + 1/7(4a^2 + 9a + 2)xy + 1/7(-4a^2 + 8a -3)y = x^3 + 1/7(-4a^2 + 8a - 3)x^2$ & 0\\ \hline
$(-a^2 + a + 2 , 2a^2 - a - 5)$ & $y^2 + 1/7(5a^2 - 8a + 6)xy + 1/7(7a^2 - 15a +5)y = x^3 + 1/7(7a^2 - 15a + 5)x^2$ & 0\\ \hline \hline
\multicolumn{3}{|l|}{Number field generated by a root of the polynomial $x^4-2x^3+x^2-2x+1$} \\ \hline

  Point from $X_1(14)$ &Corresponding Elliptic Curve $E$ with torsion $\mathbb{Z} \slash 14\mathbb{Z}$ & Rank of $E$ \\ 
  \hline
$(a^3 - a^2 + a - 1 , a^3 - 1)$ & $y^2 + 1/7(2a^3 + a^2 + a + 9)xy + 1/7(-3a^3 +3a^2 - 2a + 5)y = x^3 + 1/7(-3a^3 + 3a^2 - 2a + 5)x^2$& 0 \\ \hline
$(-a^3 + 2a^2 ,a^3 - 3a^2 + 3a - 2 )$ & $ y^2 + 1/7(-a^3 + 10a^2 - 18a + 13)xy + 1/7(3a^3- 7a^2 + 4a - 2)y = x^3 + 1/7(3a^3 - 7a^2 + 4a - 2)x^2$& 0 \\ \hline \hline

\multicolumn{3}{|l|}{Number field generated by a root of the polynomial $x^6-3x^5+4x^4-3x^3-8x^2+9x+27$} \\ \hline

  Point from $X_1(14)$ &Corresponding Elliptic Curve $E$ with torsion $\mathbb{Z} \slash 14\mathbb{Z}$& Rank of $E$  \\ 
  \hline
$(1/3(a^2 - a) , 1/225(-2a^5 + 5a^4 - 43a^3 +22a^2 + 27a -117))$ & $y^2 + 1/1050(-66a^5 + 115a^4 + 156a^3 - 199a^2 -384a + 1239)xy + 1/7350(173a^5 + 380a^4 - 1943a^3 - 1178a^2 + 3477a+ 5508)y = x^3 + 1/7350(173a^5 + 380a^4 - 1943a^3 - 1178a^2 + 3477a +5508)x^2$& 0 \\ \hline
$(1/3(a^2 - a) , 1/225(2a^5 - 5a^4 + 43a^3 - 97a^2 + 48a -108))$ & $ y^2 + 1/1050(66a^5 - 215a^4 + 44a^3 + 299a^2 +184a + 861)xy + 1/7350(-173a^5 + 1245a^4 - 1307a^3 - 2997a^2 +2323a + 6417)y = x^3 + 1/7350(-173a^5 + 1245a^4 - 1307a^3 - 2997a^2 +
2323a + 6417)x^2$&0\\ \hline
$(1/150(-11a^5 + 40a^4 - 49a^3 - 4a^2 + 261a -456) , 1/150(44a^5 - 185a^4 + 396a^3 - 559a^2 + 106a + 549) )$ & $ y^2 + 1/1575(-88a^5 + 420a^4 - 967a^3 + 1443a^2 -787a + 852)xy + 1/66150(433a^5 - 1245a^4 + 2497a^3 - 3963a^2 -1583a + 12393)y = x^3 + 1/66150(433a^5 - 1245a^4 + 2497a^3 - 3963a^2- 1583a + 12393)x^2$&0 \\ \hline
$1/150(11a^5 - 15a^4 - a^3 - 21a^2 - 211a - 219) , 1/150(-44a^5 + 35a^4- 96a^3 - 41a^2 + 344a +351))$ & $y^2 + 1/1575(88a^5 - 20a^4 + 167a^3 + 182a^2 -438a + 873)xy + 1/66150(-433a^5 + 920a^4 - 1847a^3 + 388a^2 + 4833a+ 8532)y = x^3 + 1/66150(-433a^5 + 920a^4 - 1847a^3 + 388a^2 + 4833a
+ 8532)x^2$& 0 \\ \hline \hline

\multicolumn{3}{|l|}{Number field generated by a root of the polynomial $x^6-3x^5+4x^4-x^3-2x^2+x+1$} \\ \hline

  Point from $X_1(14)$ &Corresponding Elliptic Curve $E$ with torsion $\mathbb{Z} \slash 14\mathbb{Z}$& Rank of $E$  \\ 
  \hline
$(-a^5 + 3a^4 - 4a^3 + 2a^2 , a^5 - 4a^4 + 7a^3 - 5a^2 + 1)$ & $y^2 + 1/7(-2a^4 + 10a^3 - 12a^2 + 2a + 12)xy +1/49(-3a^5 + 19a^4 - 10a^3 - 22a^2 + 21a + 11)y = x^3 + 1/49(-3a^5 +19a^4 - 10a^3 - 22a^2 + 21a + 11)x^2$ (Torsion is  $\mathbb{Z} \slash 2\mathbb{Z} \times \mathbb{Z} \slash 14\mathbb{Z}$ )& 0 \\ \hline
$(-a^5 + a^4 - a^3 - a^2 - a , -a^5 + 4a^4 - 3a^3 + 5a + 2)$ & $ y^2 + 1/21(17a^5 - 26a^4 + 24a^3 + 22a^2 - 7a +17)xy + 1/21(a^5 - 5a^4 + 6a^3 - a^2 - 6a)y = x^3 + 1/21(a^5 - 5a^4+ 6a^3 - a^2 - 6a)x^2$ & 0\\ \hline
$(a^5 - 4a^4 + 7a^3 - 6a^2 + a + 1 , -a^5 + 3a^4 - 3a^3 - a^2 + 5a - 3)$ & $ y^2 + 1/7(2a^5 - 11a^4 + 24a^3 - 22a^2 + 3a +20)xy + 1/49(9a^5 - 15a^4 + 30a^3 - 11a^2 + 14a + 23)y = x^3 +1/49(9a^5 - 15a^4 + 30a^3 - 11a^2 + 14a + 23)x^2$& 0 \\ \hline
$(a^4 - 2a^3 + 2a^2 + a - 1 , a^5 - 2a^4 + 2a^3 + a^2 - 2a)$ & $y^2 + 1/7(10a^5 - 33a^4 + 52a^3 - 27a^2 - 14a +31)xy + 1/7(-2a^5 + 5a^4 - a^3 -14a^2 + 24a - 12)y = x^3 +1/7(-2a^5 + 5a^4 - a^3 -14a^2 + 24a - 12)x^2$& 0 \\ \hline \hline
\multicolumn{3}{|l|}{Number field generated by a root of the polynomial $x^6-x^5+x^4-x^3+x^2-x+1$} \\ \hline

  Point from $X_1(14)$ &Corresponding Elliptic Curve $E$ with torsion $\mathbb{Z} \slash 14\mathbb{Z}$& Rank of $E$  \\ 
  \hline
$(a^5 + a^3 - a^2 - 1 , -a^5 + a^4 + 2a^2 + 1)$ & $y^2 + 1/7(-9a^5 - 2a^4 - 9a^3 - 6a + 13)xy +
1/7(3a^5 + a^4 + 3a^3 + a - 1)y = x^3 + 1/7(3a^5 + a^4 + 3a^3 + a -1)x^2$ (Torsion is  $\mathbb{Z} \slash 2\mathbb{Z} \times \mathbb{Z} \slash 14\mathbb{Z}$ )& 0 \\ \hline
$(1/4(-a^4 - a^2 + a + 1) , 1/8(a^4 + a^2 - a - 5))$ & $y^2 + 1/28(a^4 + a^2 - a + 32)xy + 1/56(a^4 + a^2- a + 6)y = x^3 + 1/56(a^4 + a^2 - a + 6)x^2$& 0 \\ \hline
$(2a^5- 4a^4 + 4a^3 - 2a^2 + 1 , 6a^5 - 10a^4 + 10a^3 - 6a^2 + 2)$ & $ y^2 + 1/7(4a^5 - 13a^4 + 13a^3 - 4a^2 + 6)xy +1/7(-4a^5 - 4a^4 + 4a^3 + 4a^2 - 7)y = x^3 + 1/7(-4a^5 - 4a^4 + 4a^3 +4a^2 - 7)x^2$& 0 \\ \hline
$(-a^4 + a^3 + a - 1 , -a^5 + a^4 + a^2 - 2a)$ & $y^2 + 1/14(13a^5 - 19a^4 + 4a^3 - 3a^2 + 16a +6)xy + 1/28(7a^5 + 2a^4 - 6a^3 - 2a^2 + a + 9)y = x^3 + 1/28(7a^5+ 2a^4 - 6a^3 - 2a^2 + a + 9)x^2$& 0 \\ \hline
$(-a^4 - a^2 + a - 1 , a^4 + a^2 - a + 2)$ & $ y^2 + 1/7(6a^4 + 6a^2 - 6a + 17)xy + 1/7(2a^4
+ 2a^2 - 2a - 2)y = x^3 + 1/7(2a^4 + 2a^2 - 2a - 2)x^2$& 0 \\ \hline
$(a^4 - a^3 + a^2 , -a^5 - a^2 + a - 1)$ & $ y^2 + 1/14(-3a^5 + 6a^4 - 16a^3 + 19a^2 - 15a +25)xy + 1/28(3a^5 - 9a^4 + 4a^3 - 2a^2 - 4a + 7)y = x^3 +
1/28(3a^5 - 9a^4 + 4a^3 - 2a^2 - 4a + 7)x^2$&0 \\ \hline

\end{tabular}}
\label{Tab3}
  
\end{center}
\end{table}}

\tiny{\begin{table}[H]
\caption{Equations of the non-isomorphic elliptic curves over $\mathbb{Q}(a)$ with torsion subgroup $(1,15)$.}
\begin{center}
\begin{tabular}{| m{5cm} |  m{12cm} | m{1cm} | } 
  \hline
  \multicolumn{3}{|l|}{Number field generated by a root of the polynomial $x^2-x+4$} \\ \hline

  Point from $X_1(15)$ &Corresponding Elliptic Curve $E$ with torsion $\mathbb{Z} \slash 15\mathbb{Z}$ & Rank of $E$  \\ 
  \hline
$((-1/4)a , 1/8(a - 4))$ & $y^2 + 1/64(7a + 69)xy + 1/2048(79a + 93)y = x^3+ 1/2048(79a + 93)x^2$  & 0\\ \hline \hline
\multicolumn{3}{|l|}{Number field generated by a root of the polynomial $x^2-x-1$} \\ \hline

  Point from $X_1(15)$ &Corresponding Elliptic Curve $E$ with torsion $\mathbb{Z} \slash 15\mathbb{Z}$ & Rank of $E$  \\ 
  \hline
$(-a + 1 , -a + 1)$ & $y^2 + (-10a - 5)xy + (-94a - 58)y = x^3 + (-94a- 58)x^2$ &0 \\ \hline \hline
\multicolumn{3}{|l|}{Number field generated by a root of the polynomial $x^4-x^3+2x^2+x+1$} \\ \hline

  Point from $X_1(15)$ &Corresponding Elliptic Curve $E$ with torsion $\mathbb{Z} \slash 15\mathbb{Z}$  & Rank of $E$ \\ 
  \hline
$(1/4(-a^3 + a^2 - 3a - 1) , 1/8(a^3 - a^2 + 3a -3))$ & $y^2 + 1/64(7a^3 - 7a^2 + 21a + 76)xy +
1/2048(79a^3 - 79a^2 + 237a + 172)y = x^3 + 1/2048(79a^3 - 79a^2 + 237a+ 172)x^2$ & 0\\ \hline
$(1/2(-a^3 - 1) , 1/2(-a^3 - 1))$ & $ y^2 + (-5a^3 - 20)xy + (-47a^3 - 199)y = x^3 +(-47a^3 - 199)x^2 $  & 0\\ \hline \hline
\multicolumn{3}{|l|}{Number field generated by a root of the polynomial $x^4-x^3-4x^2+4x+1$} \\ \hline

  Point from $X_1(15)$ &Corresponding Elliptic Curve $E$ with torsion $\mathbb{Z} \slash 15\mathbb{Z}$  & Rank of $E$ \\ 
  \hline
$(-a^3 + 3a , -a^3 + 3a)$ & $y^2 + (-10a^3 + 30a - 15)xy + (-94a^3 + 282a -152)y = x^3 + (-94a^3 + 282a - 152)x^2$ & 0 \\ \hline
$(2a^3 + 2a^2 - 5a - 2 , -7a^3 - 6a^2 + 17a + 4 )$ & $ y^2 + (-48a^3 + 59a^2 + 180a - 232)xy + (8418a^3- 10176a^2 - 31544a + 40260)y = x^3 + (8418a^3 - 10176a^2 - 31544a +40260)x^2$ & 0 \\ \hline \hline
\multicolumn{3}{|l|}{Number field generated by a root of the polynomial $x^4-x^3+x^2-x+1$} \\ \hline

  Point from $X_1(15)$ &Corresponding Elliptic Curve $E$ with torsion $\mathbb{Z} \slash 15\mathbb{Z}$  & Rank of $E$ \\ 
  \hline
$(a - 1 , -a^3 + a^2 - a)$ & $y^2 + (-2a^3 + 5a^2 - 5a + 3)xy + (3a^3 + 5a^2- 13a + 10)y = x^3 + (3a^3 + 5a^2 - 13a + 10)x^2$  & 0\\ \hline
$(-a^3 + a^2 + 1 , -a^3 + a^2 + 1)$ & $ y^2 + (-10a^3 + 10a^2 - 5)xy + (-94a^3 + 94a^2 - 58)y = x^3 + (-94a^3 + 94a^2 - 58)x^2$ & 0 \\ \hline

\end{tabular}
\label{Tab4}
\end{center}
\end{table}}

\tiny{\begin{table}[H]
\caption{Equations of the non-isomorphic elliptic curves over $\mathbb{Q}(a)$ with torsion subgroup $(2,10)$.}  
\begin{center}
\begin{tabular}{| m{5cm} |  m{12cm} | m{1cm} |  } 
  \hline
  \multicolumn{3}{|l|}{Number field generated by a root of the polynomial $x^4-2x^3+2$} \\ \hline
 
  Point from $X_1(2,10)$ &Corresponding Elliptic Curve $E$ with torsion $\mathbb{Z} \slash 2\mathbb{Z} \times \mathbb{Z} \slash 10\mathbb{Z}$ & Rank of $E$ \\ 
  \hline
$(-a^3 + 2a^2 - 1 , -2a + 3)$ & $y^2 = x^3 + 1/5(-6a^3 + 3a^2 + 12a + 3)x^2 +
1/5(-4a^3 + 4a^2 + 6)x$ & 0 \\ \hline \hline
\multicolumn{3}{|l|}{Number field generated by a root of the polynomial $x^6-3x^5+4x^4-3x^3-2x^2+3x+3$} \\ \hline

  Point from $X_1(2,10)$ & Corresponding Elliptic Curve $E$ with torsion $\mathbb{Z} \slash 2\mathbb{Z} \times \mathbb{Z} \slash 10\mathbb{Z}$ & Rank of $E$  \\ 
  \hline
$(1/3(-a^4 + 2a^3 - a^2) , 1/81(-16a^5 + 40a^4 - 44a^3 +26a^2 - 36a + 15))$ & $y^2 = x^3 + 1/10(-7a^4 + 14a^3 + 10a^2 - 17a -8)x^2 + 1/20(77a^4 - 154a^3 - 43a^2 + 120a + 87)x$ & 0 \\ \hline

\end{tabular}
\label{Tab5}
\end{center}
\end{table}}

\tiny{\begin{table}[H]
\caption{Equations of the non-isomorphic elliptic curves over $\mathbb{Q}(a)$ with torsion subgroup $(3,9)$.}  
\begin{center}
\begin{tabular}{| m{5cm} |  m{12cm} | m{1cm} |  } 
  \hline
  \multicolumn{3}{|l|}{Number field generated by a root of the polynomial $x^6-3x^5+5x^3-3x+1$} \\ \hline

  Point from $X_1(3,9)$ &Corresponding Elliptic Curve $E$ with torsion $\mathbb{Z} \slash 3\mathbb{Z} \times \mathbb{Z} \slash 9\mathbb{Z}$ & Rank of $E$ \\ 
  \hline
$(a^5 - 2a^4 - 2a^3 + 3a^2 + 2a - 1 , a^5 - a^4 - 2a^3 + a^2+ 2a - 1)$ & $y^2 + (-9a^5 + 32a^4 - 15a^3 - 50a^2 + 49a -11)xy + (-179a^5 + 634a^4 - 291a^3 - 992a^2 + 955a - 234)y = x^3 +(-179a^5 + 634a^4 - 291a^3 - 992a^2 + 955a - 234)x^2$ & 0\\ \hline \hline
\multicolumn{3}{|l|}{Number field generated by a root of the polynomial $x^6-x^3+1$} \\ \hline

  Point from $X_1(3,9)$ &Corresponding Elliptic Curve $E$ with torsion $\mathbb{Z} \slash 3\mathbb{Z} \times \mathbb{Z} \slash 9\mathbb{Z}$ & Rank of $E$ \\ 
  \hline
$(a^3 - a^2 + a - 1 , -2a^5+ a^4 + a^2 + a - 2)$ & $y^2 + (2a^5 - 6a^4 + 4a^2 + 4a - 4)xy + (13a^5
- 38a^4 + 25a^2 + 25a - 33)y = x^3 + (13a^5 - 38a^4 + 25a^2 + 25a -33)x^2$ & 0\\ \hline

\end{tabular}
\label{Tab6}
\end{center}
\end{table}}

\tiny{\begin{table}[H]
\caption{Equations of the non-isomorphic elliptic curves over $\mathbb{Q}(a)$ with torsion subgroup $(6,6)$.}  
\begin{center}
\begin{tabular}{| m{5cm} |  m{12cm} | m{1cm} | } 
  \hline
  \multicolumn{3}{|l|}{Number field generated by a root of the polynomial $x^6-3x^5+5x^3-3x+1$} \\ \hline
  
  Point from $X_1(6,6)$ &Corresponding Elliptic Curve $E$ with torsion $\mathbb{Z} \slash 6\mathbb{Z} \times \mathbb{Z} \slash 6\mathbb{Z}$ & Rank of $E$  \\ 
  \hline
$(-a^4 + 3a^2 + 2a - 1 , 2a^5 - 2a^4 - 4a^3 + 2a^2 + 4a - 1)$ & $y^2 + 1/3(-2a^4 + 4a^3 - 2a + 4)xy + 1/3(-2a^4+ 4a^3 + 4a^2 - 6a + 2)y = x^3 + 1/3(-2a^4 + 4a^3 + 4a^2 - 6a +2)x^2$& 0  \\ \hline
\end{tabular}
\label{Tab7}
\end{center}
\end{table}}

\tiny{\begin{table}[H]
\caption{}

\begin{center}
\begin{itemize}
    \item $1^{st}$ column: Modular elliptic curve $X_1(m,mn)$.
    \item $2^{nd}$ column: Minimal polynomial for $a\in \overline{\mathbb{Q}}$ such that the torsion of $X_1(m,mn)$ grows in the number field $\mathbb{Q}(a)$.
    \item $3^{rd}$ column: Notation for the number field  $\mathbb{Q}(a)$.
    \item $4^{th}$ column: Growth of torsion subgroup of $X_1(m,mn)$ from $\mathbb{Q}$ to $\mathbb{Q}(a)$.
    \item $5^{th}$ column: Number of non-isomorphic elliptic curves over $\mathbb{Q}(a)$  with torsion subgroup $(m,mn)$.
    \item $6^{th}$ column: Torsion subgroup over $\mathbb{Q}(a)$ for the elliptic curve in $5^{th}$ column.
\end{itemize}
\begin{tabular}{| m{1cm}  | m{4.5cm}  |m{1.2cm}  | m{2.2cm}  | m{0.5cm} |  m{5.3cm} |  } 
  \hline 
  \multirow{3}{*}{$X_1(11)$}  & $x^3-x^2+x+1$ &$L_{11,3}^{(1,10)}$ &$(1,5)\longrightarrow(1,10)$&$1$&$(1,11)$ \\ \cline{2-6}
   & $x^5-x^4-4x^3+3x^2+3x-1$& $L_{11,5}^{(1,25)}$&$(1,5)\longrightarrow(1,25)$&$3$&$(1,11),(1,11),(1,11)$ \\ \cline{2-6}
   & $x^6-x^5+2x^4-3x^3+2x^2-x+1$&$L_{11,6}^{(2,10)}$ &$(1,5)\longrightarrow(2,10)$&$3$&$(1,11),(1,11),(1,11)$ \\ \hline \hline
   
  \multirow{6}{*}{$X_1(14)$} & $x^2-x+2$ &$L_{14,2}^{(2,6)}$&$(1,6)\longrightarrow(2,6)$&$2$&$(1,14),(1,14)$ \\ \cline{2-6}
   & $x^3-x^2-2x+1$ &$L_{14,3}^{(1,18)}$&$(1,6)\longrightarrow(1,18)$&$2$&$(1,14),(1,14)$ \\ \cline{2-6}
   & $x^4-2x^3+x^2-2x+1$ &$L_{14,4}^{(1,12)}$&$(1,6)\longrightarrow(1,12)$&$2$&$(1,14),(1,14)$ \\ \cline{2-6}
   & $x^6-3x^5+4x^4-3x^3-8x^2+9x+27$ &$L_{14,6}^{(3,6)}$&$(1,6)\longrightarrow(3,6)$&$4^{\infty}$&$(1,14),(1,14),(1,14),(1,14)$ \\ \cline{2-6}
   & $x^6-3x^5+4x^4-x^3-2x^2+x+1$& $L_{14,6}^{(1,18)}$ &$(1,6)\longrightarrow(1,18)$&$4$&$(2,14),(1,14),(1,14),(1,14)$ \\ \cline{2-6}
   & $x^6-x^5+x^4-x^3+x^2-x+1$& $L_{14,6}^{(2,18)}$&$(1,6)\longrightarrow(2,18)$&$6$&$(2,14),(1,14),(1,14),(1,14),(1,14),(1,14)$ \\ \hline \hline
   
  \multirow{6}{*}{$X_1(15)$} & $x^2-x+4$& $L_{15,2}^{(2,4)}$&$(1,4)\longrightarrow(2,4)$&$1$&$(1,15)$ \\ \cline{2-6}
   & $x^2-x-1$ & $L_{15,2}^{(1,8)}$&$(1,4)\longrightarrow(1,8)$&$1$&$(1,15)$ \\ \cline{2-6}
   & $x^2-x+1$ &$L_{15,2}^{(1,8)'}$&$(1,4)\longrightarrow(1,8)$&$0$&$-$ \\ \cline{2-6}
   & $x^4-x^3+2x^2+x+1$ &$L_{15,4}^{(2,8)}$&$(1,4)\longrightarrow(2,8)$&$2$&$(1,15), (1,15)$ \\ \cline{2-6}
   & $x^4-x^3-4x^2+4x+1$ &$L_{15,4}^{(1,16)}$&$(1,4)\longrightarrow(1,16)$&$2$&$(1,15), (1,15)$ \\ \cline{2-6}
   & $x^4-x^3+x^2-x+1$ &$L_{15,4}^{(1,16)'}$&$(1,4)\longrightarrow(1,16)$&$2$&$(1,15), (1,15)$ \\ \hline \hline
   
  \multirow{3}{*}{$X_1(2,10)$} & $x^2-x-1$ &$L_{(2,10),2}^{(2,6)}$&$(1,6)\longrightarrow(2,6)$&$0$&$-$ \\ \cline{2-6}
   & $x^4-2x^3+2$ &$L_{(2,10),4}^{(1,12)}$&$(1,6)\longrightarrow(1,12)$&$1$&$(2,10)$ \\ \cline{2-6}
   & $x^6-3x^5+4x^4-3x^3-2x^2+3x+3$& $L_{(2,10),6}^{(3,6)}$&$(1,6)\longrightarrow(3,6)$&$1$&$(2,10)$ \\ \hline \hline
   
   \multirow{4}{*}{$X_1(2,12)$} & $x^2-x+1$ &$L_{(2,12),2}^{(2,4)}$&$(1,4)\longrightarrow(2,4)$&$0$&$-$ \\ \cline{2-6}
    & $x^2+1$ &$L_{(2,12),2}^{(1,8)'}$&$(1,4)\longrightarrow(1,8)$&$0$&$-$ \\ \cline{2-6}
    & $x^2-3$ &$L_{(2,12),2}^{(1,8)}$&$(1,4)\longrightarrow(1,8)$&$0$&$-$ \\\cline{2-6}
   & $x^4-x^2+1$ &$L_{(2,12),4}^{(2,8)}$&$(1,4)\longrightarrow(2,8)$&$0$&$-$ \\ \hline \hline
   
   \multirow{5}{*}{$X_1(3,9)$}  & $x^2-x+1$ &$L_{(3,9),2}^{(3,3)}$&$(1,3)\longrightarrow(3,3)$&$0$&$-$ \\ \cline{2-6}
  & $x^3-2$ &$L_{(3,9),3}^{(1,6)}$&$(1,3)\longrightarrow(1,6)$&$0$&$-$ \\ \cline{2-6}
    & $x^3-3x-1$ &$L_{(3,9),3}^{(1,9)}$&$(1,3)\longrightarrow(1,9)$&$0$&$-$ \\ \cline{2-6}
    & $x^6-3x^5+5x^3-3x+1$& $L_{(3,9),6}^{(6,6)}$&$(1,3)\longrightarrow(6,6)$&$1$&$(3,9)$ \\ \cline{2-6}
    & $x^6-x^3+1$ &$L_{(3,9),6}^{(3,9)}$&$(1,3)\longrightarrow(3,9)$&$1$&$(3,9)$ \\ \hline \hline
    
   \multirow{3}{*}{$X_1(4,8)$} & $x^2+1$ &$L_{(4,8),2}^{(2,4)'}$&$(2,2)\longrightarrow(2,4)$&$0$&$-$ \\ \cline{2-6}
     & $x^2-2$ &$L_{(4,8),2}^{(2,4)}$&$(2,2)\longrightarrow(2,4)$&$0$&$-$ \\ \cline{2-6}
     & $x^4+1$ &$L_{(4,8),4}^{(4,4)}$&$(2,2)\longrightarrow(4,4)$&$0$&$-$ \\ \hline \hline
     
    \multirow{3}{*}{$X_1(6,6)$} & $x^2-x+1$ &$L_{(6,6),2}^{(2,6)}$&$(1,6)\longrightarrow(2,6)$&$0$&$-$ \\ \cline{2-6}
     & $x^4-2x^3-2x+1$ &$L_{(6,6),4}^{(1,12)}$&$(1,6)\longrightarrow(1,12)$&$0$&$-$ \\\cline{2-6}
     & $x^6-3x^5+5x^3-3x+1$ &$L_{(6,6),6}^{(6,6)}$&$(1,6)\longrightarrow(6,6)$&$1$&$(6,6)$ \\ \hline
  \end{tabular}
  \label{Tab1}
\end{center}
\end{table}}

\nocite{*}
\bibliographystyle{plain}
\bibliography{main.bib}

\end{document}